\documentclass[12pt]{article}
\usepackage{theorem, amsfonts,amsmath,epsfig}

\nonstopmode
\def\CO{{\cal O}}
\def\CA{{\cal A}}
\def\CB{{\cal B}}
\def\CC{{\cal C}}
\def\CD{{\cal D}}

\def\CF{{\cal F}}

\def\CH{{\cal H}}

\def\CM{{\cal M}}

\def\CR{{\cal R}}

\def\CS{{\cal S}}

\def\CX{{\cal X}}

\newcommand\COMP{\hbox{C\kern -.58em {\raise .54ex
\hbox{$\scriptscriptstyle |$}}\kern-.55em {\raise .53ex
\hbox{$\scriptscriptstyle |$}} }}

\newcommand\MM{\hbox{I\kern-.2em\hbox{M}}}
\newcommand\NN{\hbox{I\kern-.2em\hbox{N}}}
\newcommand\RR{\mathbb{R}}
\newcommand\sRR{{\it \hbox{I\kern-.2em\hbox{R}}}}
\newcommand\QQ{\hbox{I\kern-.53em\hbox{Q}}}
\newcommand\PP{\hbox{I\kern -.2em\hbox{P}}}
\newcommand\EE{\hbox{I\kern-.2em\hbox{E}}}
\newcommand\ZZ{{{\rm Z}\kern-.28em{\rm Z}}}
\newcommand\II{{{\rm I}\kern-.28em{\rm I}}}

\DeclareMathOperator*{\argmin}{arg\,min}

\newcommand\qed{\hfill$\sqcap\kern-8.0pt\hbox{$\sqcup$}$}
\newcommand\be{\begin{equation}}
\newcommand\ee{\end{equation}}

\newtheorem{theorem}{Theorem}[section]
\newtheorem{proposition}[theorem]{Proposition}
\newtheorem{lemma}[theorem]{Lemma}

\newtheorem{assumption}{Assumption}

\newtheorem{corollary}[theorem]{Corollary}
\newtheorem{definition}[theorem]{Definition}

\title{A Monte Carlo method for exponential hedging of contingent claims}
\author{M. R. Grasselli$\ \!^*$
and T. R. Hurd\thanks{{Research supported by the Natural Sciences and
Engineering Research Council of Canada and Mathematics of Information Technology and Complex Systems, Canada}} \\
Dept. of Mathematics and Statistics\\
McMaster University\\Hamilton ON L8S 4K1}
\begin{document}
\maketitle

\begin{abstract} 
Utility based methods provide a very general theoretically consistent approach to pricing and hedging
of securities in incomplete financial markets. Solving problems in the utility based framework typically
involves dynamic programming, which in practise can be difficult to implement. This article presents a Monte Carlo
approach to optimal portfolio problems for which the dynamic programming is based on the exponential utility
function $U(x)=-\exp(-x)$. The algorithm, inspired by the Longstaff-Schwartz 
approach to
pricing American options by Monte Carlo simulation, involves learning the optimal portfolio selection strategy
on simulated Monte Carlo data. It shares with the LS framework intuitivity, simplicity and flexibility.
\end{abstract}

\newpage

\section{Introduction}

As realized in the pioneering work of Black, Scholes, Merton and others, 
financial assets in complete markets can be priced uniquely by 
construction of replicating portfolios and application of the no 
arbitrage principle. This conceptual framework forms the basis of much 
of the currently used methodology for financial engineering.  In recent 
years, however, finance practitioners have been increasingly led by 
competitive pressures to the use of much more general incomplete market 
models, such as those driven by noise with stochastic volatility, jumps 
or general L\'evy processes. In incomplete markets, matters are much more 
complicated, and the pricing and hedging of financial assets depends on 
the risk preferences of the investor.

Utility based portfolio theory provides a coherent, general and 
economically sound approach to risk--management in general financial 
models. This theory is built on the principle that market agents invest 
rationally by seeking to maximize their expected utility over some time 
period, where their utility function encodes the ``happiness'' they derive 
in holding a given level of wealth. Key works in this program are those of \cite{KarLehShrXu91}, 
\cite{KarShr98}, \cite{KraSch99}. The culmination of these results is a 
body of theory which give necessary and sufficient conditions for 
existence and uniqueness of optimal portfolios in a broad range of 
contexts.

Utility based pricing and hedging are extensions growing naturally out 
of portfolio optimization, and  much work is now in progress to place 
these methods in the broadest context, and to explore their various 
ramifications. The basic problem is that of a rational agent who seeks 
to find their optimal 
hedging portfolio when they have sold (or bought) a contingent claim. This 
framework leads to new concepts, notably the Davis price \cite{Davi97} and the 
indifference price of the contingent claim \cite{HodNeu89}.

This much more general theory is naturally applicable in areas such as 
insurance where the complete market theory appears inappropriate \cite{YouZar02}. In this context, 
the indifference price can be thought of as the reservation price of the 
claim, that is the amount the insurer should set aside to deal with its 
future liability.

Practical implemention of incomplete market models based on these new 
theoretical developments requires the development of efficient numerical 
methods. Three distinct approaches can be considered and ultimately all 
three are needed for a complete understanding of implementation issues. 
One approach is the numerical solution of general 
Hamilton--Jacobi--Bellman equations, which are the partial differential 
equations derived from stochastic control theory. A 
second approach could be broadly classified as ``state space 
discretization'', by which we mean tree and lattice based methods. A 
third broad approach can be called Monte Carlo or random simulation 
based methods. It is this third approach we attempt to realize in the 
present paper.

To our knowledge, Monte Carlo methods, although widely used for pricing derivatives \cite{BoyBroGla97},
have not been extensively used for optimal portfolio theory. Some works related to this in the 
context of complete markets
are \cite{DetGarRin03} and \cite{CviGouZap00}. Our proposed application of Monte Carlo is 
intrinsically more difficult than for example its use in the pricing of 
American style options, a problem which has only quite recently been 
efficiently implemented with the least squares algorithm of 
\cite{LonSch01}. Despite these difficulties, which we will see quite 
clearly in this paper, Monte Carlo methods have a great asset in being 
very simple and intuitive. By implementing such methods, we can 
gain key intuition and understanding which may be quite difficult to 
learn from the abstract theory.

The paper is organized as follows. Section 2 provides the reader with a 
rather detailed survey of the current theory of optimal portfolios. We 
give careful statements of the main results concerning the existence and 
uniqueness of optimal solutions for Merton's problem. We also review the 
framework of utility based hedging, introducing the key 
concepts and the basic existence/uniqueness results. The special case of exponential utility 
is discussed in some detail, 
because it has the important property that optimal solutions are 
independent of the level of wealth. This property has an important 
implication for our proposed Monte Carlo algorithm.

Section 3 focuses on the dynamics of portfolio optimization, in 
particular, the principle of dynamic programming. The concepts of 
certainty equivalent value, indifference price and the Davis price are 
introduced. The example of the geometric Brownian motion market is 
worked out in some detail. Section 4 specializes to the discrete time 
hedging framework and gives explicit formulas for dynamic programming.

The main innovation of the paper is the exponential utility algorithm  
given in section 5. It is a Monte Carlo method for learning the optimal 
trading strategy for the class of discrete time hedging problems 
introduced in section 4. This algorithm is inspired by the least-squares 
algorithm of Longstaff and Schwartz for pricing American options. Interestingly, our method 
works well only for the expopnenial utility, and no simple extension suggests itself for general
utility functions.
Section 6 describes our first application of the algorithm to hedging in 
a one-dimensional geometric Brownian motion model. We focus on this 
exactly solvable model in order to have explicit formulas with which to 
compare our Monte Carlo simulation. While the hedging strategies learned 
by the algorithm are somewhat crude, we find that the computed 
indifference prices are quite accurate. In our concluding section 7, we 
discuss the various advantages and drawbacks we see in the method.

\section{Utility based hedging for semimartingale markets}

The hedging problem is the problem of a market agent who faces a liability $B$ at a
time $T$
and must invest in the market over the period $[0,T]$ in an efficient,
rational or otherwise optimal way to reduce the risk of the liability. The randomness of the
market is represented by a filtered probability space $(\Omega,\CF,(\CF_t)_{t\in[0,T]},P)$ satisfying the 
``usual conditions'' of right continuity and completeness and we assume for simplicity that $\CF=\CF_T$. 
The discounted prices of tradeable assets 
in the market are given by the $\RR^d$--valued {\it c\`adl\`ag} semimartingale $S_t=(S^1_t,\ldots,S^d_t)$ on the 
filtration $(\CF_t)$. The liability $B$ is assumed to be an $\CF_T$--measurable random variable.

A portfolio process, or a trading strategy, is an $\RR^d$-valued predictable $S$--integrable process
$H_t=(H^1_t,\ldots,H^d_t)$, which represents the agent's asset allocations, that is, how many units
of each traded asset are held by the agent at each time $t$. The class
of such processes is denoted by $L(S)$ \cite{Prot90}. We assume that the portfolio
is {\it self-financing} (i.e. the changes in its discounted market value are solely due to the random 
changes in the prices of the
traded assets) so that the agent's discounted wealth at each time $t$ is given by the process
\[X_t=x+(H\cdot S)_t:=x+\int_0^t H_u dS_u, \qquad t\in[0,T],\]
where $x\in\RR$ is some deterministic initial wealth.    
 
To rule out strategies for which the wealth assumes arbitrarily negative values (such as ``doubling strategies''), 
we need to assume some admissibility condition on the possible portfolio
processes. Following 
\cite{HarPli81}, we say that
\begin{definition} The class $\CH$ of admissible portfolios consists of the process $H\in L(S)$ for which 
$(H\cdot S)_t$ is $P$--a.s. uniformly bounded from below. 
\end{definition}

More explicitly, $H$ is admissible if there exists a constant $k\geq 0$ (possibly depending on $H$, but neither
on $t$ nor on $\omega$) such that
\[(H\cdot S)_t(\omega) \geq -k,\]
for almost all $\omega\in \Omega$ and all $t\in [0,T]$. 

As a first consequence  of this notion of admissibility,
we have the following useful result concerning the closedness of the class of local martingales under stochastic integration
\cite[theorem 2.9]{DelSch94}:
\begin{lemma} If $S$ is a local martingale and $H$ is an admissible integrand for $S$, then $(H\cdot S)$ is
a local martingale. Consequently, $(H\cdot S)$ is a supermartingale.
\label{admiss}
\end{lemma}

Regarding martingale measures, we adopt the following definition.
\begin{definition} A probability measure $Q$ is called an absolutely continuous (resp. equivalent) 
local martingale measure for $S$ if $Q\ll P$ (resp. $Q\sim P$) and $S$ is a local martingale under $Q$.
\end{definition}

We denote the set of absolutely continuous (resp. equivalent) local martingale measures for the price process $S$ by 
$\CM^a(S)$ (resp. by $\CM^e(S)$). Observe that, due to lemma \ref{admiss}, 
a probability measure $Q\ll P$ (resp. $Q\sim P$) is an absolutely continuous (resp. equivalent) local martingale
measure if and only if $(H\cdot S)$ is a local martingale under $Q$ for any $H\in\CH$.

To ensure a viable market, free of arbitrage, we assume the technical condition ``No Free Lunch with
Vanishing Risk'' (NFLVR), which is slightly more general than ``No Arbitrage'' (NA). The reader is
referred to \cite[sections 2 and 3]{DelSch94} for the precise definition of these notions, as
well as the relations between them. 
In its most general form \cite{DelSch98}, the fundamental theorem of 
asset pricing (FTAP) asserts the equivalence between (NFLVR) and the existence of an equivalent 
$\sigma$--martingale measure for the price process $S$ \cite{DelSch98}, which might
fail to be in  $\CM^e(S)$ if we allow $S$ to have unbounded unpredictable jumps. The technicality of
using $\sigma$--martingales can be avoided, however, if we restrict ourselves to price processes $S$ which
are {\it locally bounded}. By that we mean that there exists a localizing sequence of stopping times $\{T_n\}$ 
such that, for each $n$, the stopped processes $S^{T_n}$ are bounded. In this context, we have
\cite[corollary 1.2]{DelSch94}:
\begin{theorem}[FTAP] If $S$ is a locally bounded semimartingale, then there exists an equivalent local
martingale measure $Q$ for $S$ if and only if $S$ satisfies (NFLVR).
\end{theorem}

In view of this theorem, we will henceforth assume that $S$ is locally bounded and that
\begin{assumption}[NFLVR] $\CM^e(S)\ne \emptyset$.
\label{NFLVR}
\end{assumption} 

The hedging problem can be made specific by introducing the agent's utility $U:\RR\to\RR\cup\{-\infty\}$,
a concave, strictly increasing, differentiable function. Beginning with initial capital $x\in\RR$, 
the agent then solves the optimal hedging problem
\be
\sup_{H\in\CH} E\left[U\left(x+(H\cdot S)_T - B\right)\right].
\label{optimalhedge}
\ee

If $B\equiv 0$, the optimal hedging problem reduces to Merton's optimal investment problem
\be
\sup_{H\in\CH} E\left[U\left(x+(H\cdot S)_T\right)\right].
\label{Merton}
\ee

To assert the existence and uniqueness of solutions to problems of the form (\ref{Merton}) in 
incomplete markets, one first needs to impose further technical restrictions on the class of utility functions.
In the next assumption we summarize the main properties required to hold throughout this paper. They include 
the ``reasonable asymptotic elasticity'' condition as defined in \cite[definition 1.5]{Scha01a}.

\begin{assumption}
The utility function $U:\RR\rightarrow\RR\cup\{-\infty\}$ is increasing on 
$\RR$, continuous on $\{U>-\infty\}$, differentiable and strictly concave on 
$\mbox{dom}(U)=\mbox{int}\{U>-\infty\}$,
satisfying
\begin{equation}
\lim_{x\rightarrow\infty} U^\prime(x)=0.
\end{equation}
Furthermore, we assume that one of the following cases hold.

\noindent
{\bf Case 1:} $\mbox{dom}(U)=(0,\infty)$, with 
$\displaystyle{\lim_{x\rightarrow 0} U^\prime(x)=\infty}$ and 
$\displaystyle{\limsup_{x\to\infty}} \frac{xU^\prime(x)}{U(x)}<1$.

\noindent
{\bf Case 2:} $\mbox{dom}(U)=\RR$, with $\displaystyle{\lim_{x\rightarrow -\infty} U^\prime(x)=\infty}$, 
\[\limsup_{x\to\infty} \frac{xU^\prime(x)}{U(x)}<1 \quad \mbox{and} \quad 
\liminf_{x\to-\infty} \frac{xU^\prime(x)}{U(x)}>1.\]
\label{utilass}
\end{assumption}

The central technical weaponry used to address the general solution to problem (\ref{Merton}) is convex duality,
by means of which the utility maximization problem over admissible portfolios (the ``primal problem'') is related to a 
minimization problem over a suitable domain in the set of measures on $\Omega$ (the ``dual problem''). 
The first step is to define the conjugate function $V$ as the Legendre 
transform of the function $-U(-x)$, that is
\begin{equation}
\label{conjugate}
V(y):=\sup_{x\in\RR}[U(x)-xy],\qquad y>0.
\end{equation}

It follows from well known results in convex analysis \cite{Rock97}, that the function $V$ has the 
properties listed below.

\begin{proposition} If $U$ satisfies assumption \ref{utilass}, then the 
conjugate function $V$ is
finite valued, differentiable, strictly convex on $(0,\infty)$ and 
satisfies
\begin{equation}
\lim_{y\rightarrow 0} V (y)=\lim_{x\rightarrow\infty} U(x), \qquad 
\lim_{y\rightarrow 0} V^\prime (y)=-\infty.
\end{equation}
Moreover, the behaviour of $V$ at infinity is determined by the two cases 
in assumption \ref{utilass} as follows:

\noindent
{\bf Case 1:}
$\displaystyle{\lim_{y\rightarrow\infty} V(y)=\lim_{x\rightarrow 0} U(x) 
\quad \mbox{and } \quad
\lim_{y\rightarrow\infty} V^\prime (y)=0}$.

\noindent
{\bf Case 2:}
$\displaystyle{\lim_{y\rightarrow\infty} V(y)=\infty \qquad \qquad 
\mbox{and } \quad
\lim_{y \rightarrow\infty} V^\prime (y)=\infty}$.
\label{vprop}
\end{proposition}

Both the primal and dual problems are solved over different domains depending on
which of the two cases above we are dealing with. We start with the first case, for which the present
state--of--the--art solution can be found in \cite{KraSch99}. Since in this case the utility function is
only defined for positive wealths, we need to consider the set
\begin{equation}
\CX(x)=\{X\geq 0 : X_t=x+(H\cdot S)_t, \mbox{ for some } H\in L(S), 0\leq t\leq T\}.
\end{equation} 
It is clear that $x+(H\cdot S)_t\geq 0$ implies that the portfolio $H$ must be admissible, that is
\[\CX(x)\subset \{X_t=x+(H\cdot S)_t, H\in\CH, t\in[0,T]\}\]
with a strict inclusion. Next we move from the set of processes $\CX(x)$ to the set of positive random variables
\begin{equation}
C(x)=\{g\in L^0_+(\Omega,\CF_T,P): g\leq X_T, \mbox{ for some } X\in \CX(x)\}
\label{cmerton}
\end{equation}
and observe that, since the utility function is increasing, the primal problem for case 1 written in the form
\be
\sup_{X\in\CX(x)} E\left[U(X_T)\right].
\label{Mertonx}
\ee
is equivalent to
\be
u(x)=\sup_{g\in C(x)} E\left[U(g)\right].
\label{Mertonc}
\ee

At this point, in order to exclude trivial cases, we make the following assumption.
\begin{assumption}
The value function $u$ defined in (\ref{Mertonc}) satisfies $u(x)<\infty$, for some $ x>0$.
\label{finite1}
\end{assumption}

As for the domain of the dual problem, one looks for a set with
the property of being in a ``polar relation'' with the set $C$ (the reader is refered to 
\cite{BraSch99} for the definition of the polar of a subset of $L^0_+(\Omega,\CF,P)$). 
In the mathematical finance literature \cite{DelSch94,KarLehShrXu91}, the sets
$\CM^e(S)$ and $\CM^a(S)$ were considered. One of the main technical novelties in
\cite{KraSch99} was to enlarge this domain in order to obtain a set $D$ in ``perfect'' polar relation with 
$C$ (see \cite[proposition 3.1]{KraSch99}). The set $D$ turns out to be the
convex, solid hull of $\CM^a(S)$ in $L^0_+(\Omega,\CF,P)$ (topologized by convergence in measure). Amongst
the several equivalent characterizations of the set $D$, we single out the following 
\cite{Scha01b}
\begin{eqnarray}
D&=&\left\{Y_T\in L^0_+(\Omega, \CF_T,P): \mbox {there exists a sequence } \right. \nonumber \\
&& \left. (Q_n)_{n=1}^\infty \in \CM^a(S) \mbox{  such that } Y_T\leq (a.s.)\lim_{n\to\infty} \frac{dQ_n}{dP} \right\}
\end{eqnarray}

For $y>0$, let us define $D(y)=yD$. The dual problem for case 1 can now be formulated as
\be
v(y)=\inf_{Y_T\in D(y)} E[V(Y_T)].
\label{duald}
\ee

The next theorem states the existence and uniqueness of solution for the problems (\ref{Mertonc}) and (\ref{duald})
for utilities restricted to positive wealth \cite[theorem 2.2]{KraSch99}.

\begin{theorem}  
\label{case1}
Suppose that assumptions \ref{NFLVR}, \ref{utilass} (case 1) and \ref{finite1} are satisfied. 
Then, for any $x\in\mbox{dom}(U)$ and $y>0$, the problems
\be
u(x)=\sup_{X_T\in C(x)} E[U(X_T)], \qquad v(y)=\inf_{Y_T\in D(y)}E[V(Y_T)]
\ee
have unique optimizers $\widehat{X}_T(x)\in C(x)$ and $\widehat{Y}_T(y)\in D(y)$ satisfying
\be
U^\prime(\widehat{X}_T(x))=\widehat{Y}_T(y),
\label{rel1}
\ee
where $x$ and $y$ are related by $u^\prime(x)=y$. 
\end{theorem}

We note that this theorem and its proof apply unchanged if we modify case 1 of
assumption \ref{utilass} to allow for utilities defined on an interval of the form $(a,\infty)$, for any $a\in\RR$,
provided we impose that $\displaystyle{\lim_{x\rightarrow a}} U^\prime (x)=\infty$. Observe also that the optimizer
$\widehat{X}_T(x)$ can be uniquely expressed as 
\[\widehat{X}_T(x)=x+(\widehat{H}(x)\cdot S)_T,\]
for $\widehat{H}(x)\in\CH$, whereas the 
optimizer $\widehat{Y}_T(y)$, 
even for the cases where $\widehat{Y}_T(y)/y$ is {\it not} the density of an absolutely continuous
martingale measure (by having its total $P$--mass strictly less than 1), 
can be arbitrarily approximated by elements in $\CM^a(S)$ (in the sense of almost sure convergence).

For utility functions defined on the entire real line the problem is more 
involved, due to the fact that the class of admissible portfolios as in definition \ref{admiss} 
turns out to be too narrow to contain the optimal solution. One approach is to start with the dual problem, 
for which \cite{BelFri02} shows that an optimal solution always exists (under very general
conditions). Then the set of allowed portfolios can be characterized in terms of it. This opens up a plethora of 
definitions of ``allowed'' portfolios. The reader interested in this line of thought is referred to 
\cite{DGRSSS02,KabStr02}, where the exponential utility is addressed, and to \cite{Scha02}
for more general utility functions. 

A more direct idea is to concentrate on random variables which do not necessarily arise as  
terminal values of wealth processes for any portfolios, but which can be arbitrarily approximated by such objects. 
Different such domains of optimization over random variables have been proposed \cite{Scha01a,Frit02}, 
the difference being the kind of topology (convergence)
adopted to describe the approximation mentioned above. In what follows, we adopt the approach proposed in 
\cite{Scha01a}, and specialize later on to the case of exponential utility $U(x)=-\frac{e^{-\gamma x}}{\gamma}$,
where sharper results can be quoted.

We denote by $C_U^b(x)$ the class of random variables which have integrable utility and can be dominated by the 
terminal wealth of admissible portfolios, that is,  
\begin{eqnarray}
C^b_U(x)&=&\{g\in L^0(\Omega,\CF_T,P): g\leq x+(H\cdot S)_T \nonumber \\
&& \mbox{ for some } H\in\CH \mbox{ and } U(g)\in L^1(\Omega,\CF_T,P) \}.
\label{cbMerton}
\end{eqnarray}
Next consider the closure of the set $\{U(g): g\in C^b_U(x)\}$ in the topology of $L^1(\Omega,\CF_T,P)$. Putting
$U(\infty):=\displaystyle{\lim_{x\to\infty}} U(x)$, we see that the
the utility function is a bijection between $\RR$ and $\RR$, if $U(\infty)=\infty$,
and a bijection between $\RR\cup\{\infty\}$ and $(-\infty,U(\infty)]$ otherwise. We can therefore write a general element
in this closure as $U(f)$ for some $f\in L^0(\Omega,\CF_T,P;\RR\cup\{\infty\})$. The set of such random variables
is denoted by $C_U(x)$, that is,
\begin{eqnarray}
C_U(x)&=&\left\{f\in L^0(\Omega,\CF_T,P;\RR\cup\{\infty\}): U(f) \mbox{ is in the }\right.\nonumber \\
&& \left. L^1(P)\mbox{-closure of } \{U(g): g\in C^b_U(x)\}\right\}.
\end{eqnarray}

The primal optimization problem then becomes
\be
u(x)=\sup_{f\in C_U(x)} E[U(f)].
\label{Mertonwhole}
\ee

As in case 1, to exclude trivial cases we make the following assumption.
\begin{assumption}
The value function $u$ defined in (\ref{Mertonwhole}) satisfies $u(x)< U(\infty)$, for some $x\in \RR$.
\label{finite2}
\end{assumption}
 
Complicated as the domain $C_U(x)$ might seem, the good news is that in this setting
the optimization domain for the dual problem is simply $\CM^a(S)$, as opposed to the enlarged set $D$ of case 1.
That is, the dual problem is now
\be
v(y)=\inf_{Q\in\CM^a(S)} E\left[V\left(y\frac{dQ}{dP}\right)\right].
\label{dualwhole}
\ee

We can now state a theorem for case 2 \cite[theorem 2.2]{Scha01a}.

\begin{theorem}   
\label{case2}
Suppose that assumptions \ref{NFLVR}, \ref{utilass} (case 2) and \ref{finite2} are satisfied. 
Then:
\begin{enumerate}
\item For any $x\in\RR$ and $y>0$, the problems
\be
u(x)=\sup_{f\in C_U(x)} E[U(f)], \qquad v(y)=\inf_{Q\in M^a(S)} E\left[V\left(y\frac{dQ}{dP}\right)\right]
\ee
have unique optimizers $\widehat{f}(x)\in C_U(x)$ and $\widehat{Q}(y)\in \CM^a(S)$ satisfying
\be
U^\prime(\widehat{f}(x))=y\frac{d\widehat{Q}(y)}{dP},
\label{rel2}
\ee
where $x$ and $y$ are related by $u^\prime(x)=y$.
\item If it occurs that $\widehat{Q}(y)\in\CM^e(S)$, then $\widehat{f}(x)$ equals the terminal value $\widehat{X}_T(x)$
of a uniformly integrable $\widehat{Q}(y)$-martingale of the form 
\[\widehat{X}_t (x)=x+(\widehat{H}(x)\cdot S)_t, \quad t\in[0,T],\]
for some $\widehat{H}(x)\in L(S)$.  
\end{enumerate}
\end{theorem}

Observe that the optimizer $\widehat{f}(x)\in C_U(x)$ does not need to be the terminal 
wealth of any portfolio. However, by construction of the set $C_U(x)$, its utility can be arbitrarily approximated 
by the utility of terminal wealth of admissible portfolios. As for the optimizer of the dual problem, 
recall from proposition \ref{vprop} that $\displaystyle{\lim_{y\to 0}}V(y)=U(\infty)$. Therefore for all cases 
when $U(\infty)=\infty$, the minimizer must satisfy $\frac{d\widehat{Q}(y)}{dP}>0$ almost surely, 
implying that $\widehat{Q}(y)\in\CM^e(S)$ and item 2 holds. In such cases, $\widehat{f}(x)$ itself can be 
achieved by trading according to a portfolio 
$\widehat{H}\in L(S)$. Although $\widehat{H}$ might not be in $\CH$, the wealth process generated by it, being a 
uniformly integrable $\widehat{Q}(y)$ martingale, certainly does not arise from a ``doubling strategy'', 
so that $\widehat{H}$ can be considered {\it a posteriori} to be an 
``allowed'' portfolio. Turning this argument around was the starting point of the aforementioned approaches to extend 
the domain of the primal problem to include such portfolios \cite{DGRSSS02,KabStr02,Scha02}

But the minimizer $\widehat{Q}(y)$ is also equivalent to $P$ in other cases, 
and it is here that we specialize to an exponential utility of the form 
$U(x)=-\frac{e^{-\gamma x}}{\gamma}$, $\gamma>0$. Observe that for this utility we have
\[\limsup_{x\to\infty} \frac{xU^\prime(x)}{U(x)}=-\infty <1\]
and
\[\liminf_{x\to -\infty} \frac{xU^\prime(x)}{U(x)}=\infty > 1,\]
so that it satisfies all the conditions for case 2 of assumption \ref{utilass}.
Observe further that 
its dual function is
\[V(y)=\frac{y}{\gamma}(\log y - 1),\]
so that the dual problem (\ref{dualwhole}) is equivalent to the problem of finding a measure in $\CM^a(S)$ with
minimal relative entropy with respect to $P$, that is,
\begin{equation}
\inf_{Q\in\CM^a(S)} E\left[\frac{dQ}{dP}\log\left(\frac{dQ}{dP}\right)\right].
\end{equation} 
It follows from \cite{Csis75} that the minimizer of this problem (which incidentally is independent of $y$)
will be an equivalent local martingale measure provided there
exists at least one measure in $\CM^e(S)$ with finite relative entropy, allowing us to use item 2 of theorem 
\ref{case2}.

\begin{corollary} Let $U(x)=-\frac{e^{-\gamma x}}{\gamma}$, $\gamma>0$, and suppose that assumptions
\ref{NFLVR} and \ref{finite2} hold. If in addition we have that  
\be
E\left[\frac{dQ}{dP}\log\left(\frac{dQ}{dP}\right)\right]<\infty,
\ee
for some $Q\in\CM^e(S)$, then the minimizer $\widehat{Q}(y)$ of theorem \ref{case2} is the equivalent 
local martingale measure $\widehat{Q}$, independent of $y>0$, which minimizes the relative entropy with respect to
P among all absolutely continuous martingale measures. Therefore $\widehat{f}(x)$ equals the terminal 
value $\widehat{X}_T(x)$
of a uniformly integrable $\widehat{Q}$-martingale of the form 
\[\widehat{X}_t (x)=x+(\widehat{H}(x)\cdot S)_t,\]
for some $\widehat{H}(x)\in L(S)$.  
\label{strongexp}
\end{corollary}

We now move to the subject of solving the hedging problem (\ref{optimalhedge}). Once more the solutions will
take place in different domains and involve different techniques depending on whether our utility function
falls into case 1 or case 2 of assumption \ref{utilass}. In either case, we are going to assume that the random 
claim that we want to hedge is a bounded random variable.

\begin{assumption} $B \in L^{\infty}(\Omega,\CF_T,P)$.
\label{boundedclaim}
\end{assumption}

We start with the first case, which was solved in \cite{CviSchWan01}. Observe that to account for the
presence of a random claim at time $T$, it is not enough to consider {\it positive} random variables which are
dominated by terminal values of admissible portfolios, as was done in (\ref{cmerton}). We therefore consider the set
\begin{equation}
\CC(x)=\{g\in L^0(\Omega,\CF_T,P):g\leq x+(H\cdot S)_T, \mbox{ for some } H\in\CH\}.
\label{chedge}
\end{equation}
The primal problem now becomes
\begin{equation}
u(x)=\sup_{g\in \CC(x)} E[U(g-B)],
\label{hedgec}
\end{equation}
where it is understood that $U(x)=-\infty$ whenever $x\leq 0$.

As in the previous cases, we assume the following.
\begin{assumption} The value function $u$ defined in (\ref{hedgec}) satisfies $|u(x)|<\infty$ for some $x>\|B\|_\infty$.
\label{finite3}
\end{assumption}

Recall that the crucial point in the proof of theorem \ref{case1} was the use of the polar relation
between the sets $C$ and $D$ as subsets of $L^0_+(\Omega,\CF_T,P)$, for which a version of the bipolar theorem can
be used \cite{BraSch99}. In the absence of such results for subsets of $L^0(\Omega,\CF_T,P)$ as a whole,
we are led to consider an appropriate subset of $L^\infty(P)$, namely 
\begin{equation}
\CC=\CC(0)\cap L^\infty(\Omega,\CF_T,P).
\label{cchedge}
\end{equation}

Accordingly, to obtain a perfect polar relation, we need to modify the definition for the domain of the
dual problem. The natural space to define the polar of a subset of $L^\infty$ is its topological dual 
$(L^\infty)^*$. We therefore define
\begin{equation}
\CD=\{Q\in (L^\infty)^* : \|Q\|=1 \mbox{ and } Q(g)\leq 0 \mbox{ for all } g\in\CC\}.
\end{equation}
To obtain a more concrete characterization of this set, notice that $\CD\in (L^\infty)^*_+$ (since $\CC$ contains
all the {\it negative} bounded random variables). The good news about the set $(L^\infty)^*_+$ is that it can
be identified with the set af all nonnegative {\it finitely} additive bounded set functions on $\CF_T$ which
vanish on the $P$--null sets. Moreover, any such function $Q\in (L^\infty)^*_+$ 
can be uniquely decomposed into its regular part $Q^r$ and its singular part $Q^s$ as follows
\[Q=Q^r+Q^s,\]
where $Q^r\geq 0$ is {\it countably} additive and $Q^s\geq 0$ is {\it purely finitely} additive. Naturally, 
$Q^r$ corresponds to a measure which is absolutely continuous with respect to $P$ and whose
Radon--Nikodym derivative is denoted by $\frac{dQ^r}{dP}$. We now look at the subset of regular elements in $\CD$, namely
\begin{equation}
\CD^r=\{Q\in\CD:Q^s=0\}=\CD\cap L^1(\Omega,\CF_T,P).
\end{equation}
Since all elements in $\CD$ have unit norm, it follows that $\CD^r$ consists of {\it probability} measures which are 
absolutely continuous with respect to $P$. In fact, since we are assuming that
the processes $S$ are locally bounded, it can be shown that $\CD^r$ is nothing but our familiar $\CM^a (S)$, 
the set of absolutely continuous local martingale measures for $S$ \cite[lemma 1.1 (b)]{BelFri02}. The dual
problem in this case is
\begin{equation}
v(y)=\inf_{Q\in\CD}\left\{E\left[V\left(y\frac{dQ^r}{dP}\right)-y\frac{dQ^r}{dP}B\right]-yQ^s(B)\right\},
\label{dualhedge1}
\end{equation}
where one should notice that the domain of optimization is the entire $\CD$, with the dual function $V$ 
contributing to it only through its regular subset $\CD^r$, whereas the dependence on the claim $B$ is manifested 
on both its regular and singular parts. In this respect, it is worth mentioning that our old set $D$ (for Merton's problem)
can also be characterized as the regular part of the weak-star closure of the convex
solid hull of $\CM^a(S)$ in $(L^\infty)^*$ (whose elements can have total $P$--mass strictly less than 1).
From this perspective, it becomes clear that the extra care necessary to treat the hedging problem in this case comes
from dealing with both the regular and singular parts of elements in the domain of the dual problem. The main result
in this case is \cite[theorem 3.1]{CviSchWan01}   

\begin{theorem}  
\label{hedge1}
Suppose that assumptions \ref{NFLVR}, \ref{utilass} (case 1), \ref{boundedclaim} and  \ref{finite3} are satisfied. 
Let $x_0=\displaystyle{\sup_{Q\in\CD}}Q(B)$. Then, for any $y>0$, the dual problem  
\be
v(y)=\inf_{Q\in\CD}\left\{E\left[V\left(y\frac{dQ^r}{dP}\right)-y\frac{dQ^r}{dP}B\right]-yQ^s(B)\right\}
\ee
has a unique (up to singular part) optimizer $\widehat{Q}(y)\in \CD$ and, for any $x>x_0$, the primal problem
\be
u(x)=\sup_{X_T\in\CC(x)} E[U(X_T-B)]
\ee
has unique optimizer $\widehat{X}_T(x)\in \CC(x)$  satisfying
\be
U^\prime(\widehat{X}_T(x)-B)=y\frac{d\widehat{Q}^r(y)}{dP},
\ee
where $x$ and $y$ are related by $u^\prime(x)=y$. 
\end{theorem}

Regarding the second case of assumption \ref{utilass}, the optimal hedging problem has been solved in \cite{DGRSSS02}
for the exponential utility and claims $B$ satisfying a boundedness conditions weaker than assumption \ref{boundedclaim}. 
In \cite{Owen02}, the problem was solved for general utility functions with reasonable asymptotic elasticity (which
include the exponential) but bounded claims (although some remarks are offered on how to extend the result to 
possibly unbounded ones). We describe here the solution of \cite{Owen02}, since it follows the same techniques 
of \cite{Scha01a} and \cite{CviSchWan01}, for
which we have already developed most of the notation. In the presence of a claim satisfying assumption
\ref{boundedclaim}, the analogue of the set $C^b_U$ defined in (\ref{cbMerton}) is
\begin{eqnarray}
\CC^b_U(x)&=&\{g\in L^0(\Omega,\CF_T,P): g\leq x+(H\cdot S)_T - B\nonumber \\
&& \mbox{ for some } H\in\CH \mbox{ and } U(g)\in L^1(\Omega,\CF_T,P) \}.
\label{cbhedge}
\end{eqnarray}
Similarly, we replace the set $C_U(x)$ by 
\begin{eqnarray}
\CC_U(x)&=&\left\{f\in L^0(\Omega,\CF_T,P;\RR\cup\{\infty\}): U(f-B) \mbox{ is in the }\right.\nonumber \\
&& \left. L^1(P)\mbox{-closure of } \{U(g): g\in C^b_U(x)\}\right\}.
\end{eqnarray}
The interpretation of this set is the same as before, only this time we have to account for the random claim $B$. 
Namely, it 
consists of random variables which, after subtracting the claim $B$, have a utility that can be arbitrarily 
approximated by the utility of terminal wealth of admissible portfolios less the claim $B$.

Our modified primal problem now reads
\be
u(x)=\sup_{f\in \CC_U(x)} E[U(f-B)],
\label{hedgewhole}
\ee
for which we assume the following.
\begin{assumption}
The value function $u$ defined in (\ref{hedgewhole}) satisfies $u(x)< U(\infty)$, for some $x\in \RR$.
\label{finite4}
\end{assumption}

As with the case of no claim, when we pass to utilities defined on the entire real line
the domain of the dual problem becomes simpler, 
being just the set $\CM^a(S)$ (as opposed to the complicated set $\CD$). In the same vein, the statement of the dual
problem is much more transparent, since it does not involve the singular measures that we encountered before. It is simply
(compare with (\ref{dualhedge1}))
\begin{equation}
v(y)=\inf_{Q\in\CM^a(S)}E\left[V\left(y\frac{dQ}{dP}\right)-y\frac{dQ}{dP}B\right].
\label{dualhedge2}
\end{equation}

The next theorem \cite[theorem 1.1]{Owen02} provides the existence and uniqueness of solutions to the 
hedging problem for utilities defined on the entire $\RR$. The remark following theorem \ref{case2}
about the optimal measure $\widehat{Q}(y)$ being actually equivalent to $P$ when $U(\infty)=\infty$ applies here as well
(as can be seen from the form of the dual problem (\ref{dualhedge2})). 

\begin{theorem}   
\label{hedge2}
Suppose that assumptions \ref{NFLVR}, \ref{utilass} (case 2), \ref{boundedclaim} and \ref{finite4} are satisfied. 
Then:
\begin{enumerate}
\item For any $x\in\RR$ and $y>0$, the problems
\[u(x)=\sup_{f\in \CC_U(x)} E[U(f-B)], \quad v(y)=\inf_{Q\in M^a(S)} E\left[V\left(y\frac{dQ}{dP}\right)-
y\frac{dQ}{dP}B\right]\]
have unique optimizers $\widehat{f}(x)\in C_U(x)$ and $\widehat{Q}(y)\in \CM^a(S)$ satisfying
\be
U^\prime(\widehat{f}(x)-B)=y\frac{d\widehat{Q}(y)}{dP},
\ee
where $x$ and $y$ are related by $u^\prime(x)=y$.
\item If it occurs that $\widehat{Q}(y)\in\CM^e(S)$, then $\widehat{f}(x)$ equals the terminal value $\widehat{X}_T(x)$
of a uniformly integrable $\widehat{Q}(y)$-martingale of the form 
\[\widehat{X}_t (x)=x+(\widehat{H}(x)\cdot S)_t,\]
for some $\widehat{H}(x)\in L(S)$.  
\end{enumerate}
\end{theorem}

To assert that the optimal measure $\widehat{Q}(y)$ is actually equivalent to $P$ for the case of exponential utility,
the analogue of proposition \ref{strongexp}, we follow \cite{DGRSSS02} and
consider the change from $P$ to an equivalent probability measure $P_B$ with density
\be
\frac{dP_B}{dP}=c_Be^{\gamma B}, \mbox{ with } c_B^{-1}=E[e^{\gamma B}].
\ee
Therefore, for any $Q\ll P$, we have that
\be
E\left[\frac{dQ}{dP}\log\frac{dQ}{dP}\right]=E_{P_B}\left[\frac{dQ}{dP_B}\log\frac{dQ}{dP_B}\right]+\log c_B +
E\left[\gamma\frac{dQ}{dP}B\right].
\ee
It then follows from the boundedness of $B$ that $Q$ has finite relative entropy with respect to $P$ if and only
if it has finite relative entropy with respect to $P_B$.

Now notice that the dual problem in this case is
\begin{eqnarray}
v(y)&=&\inf_{Q\in\CM^a(S)}E\left[\frac{y}{\gamma}\frac{dQ}{dP}\left(\log\left(y\frac{dQ}{dP}\right)-1\right)
-y\frac{dQ}{dP}B\right] \nonumber \\
&=& \frac{y}{\gamma}(\log y -1)+\inf_{Q\in\CM^a(S)}E\left[\frac{dQ}{dP}\log\frac{dQ}{dP}-\gamma\frac{dQ}{dP}B\right]
\nonumber \\
&=& \frac{y}{\gamma}(\log yc_B -1)+ \frac{y}{\gamma}\inf_{Q\in\CM^a(S)}E_{P_B}
\left[\frac{dQ}{dP_B}\log\frac{dQ}{dP_B}\right],
\end{eqnarray}
from which we see that its minimizer coincides with the minimizer of the relative entropy with respect to $P_B$
over all the absolutely continuous martingale measures for $S$. But from the argument preceding corollary
\ref{strongexp}, such a minimizer is equivalent to
$P_B$ (and therefore to $P$) provided there is at least one $Q$ in $\CM^e(S)$ whose relative entropy with respect to $P_B$ 
is finite, which in turn is the same as having at least one $Q$ in $\CM^e(S)$ whose relative entropy with respect to $P$ is
finite. This suffices to prove:

\begin{corollary} Let $U(x)=-\frac{e^{-\gamma x}}{\gamma}$, $\gamma>0$, and suppose that assumptions
\ref{NFLVR}, \ref{boundedclaim} and \ref{finite4} hold. If in addition we have that  
\be
E\left[\frac{dQ}{dP}\log\left(\frac{dQ}{dP}\right)\right]<\infty,
\ee
for some $Q\in\CM^e(S)$, then the minimizer $\widehat{Q}(y)$ of theorem \ref{hedge2} is the equivalent 
local martingale measure $\widehat{Q}$, independent of $y>0$, which minimizes the relative entropy with respect to
$P_B$ among all absolutely continuous martingale measures. Therefore $\widehat{f}(x)$ equals the terminal 
value $\widehat{X}_T(x)$
of a uniformly integrable $\widehat{Q}$-martingale of the form 
\[\widehat{X}_t (x)=x+(\widehat{H}(x)\cdot S)_t,\]
for some $\widehat{H}(x)\in L(S)$.  
\end{corollary}

We end this review section with a discussion about complete markets, defined to be those for which 
there is exactly one equivalent martingale measure $Q$, that is, $\CM^e(S)$ is the singleton $\{Q\}$. 
The second fundamental theorem of asset pricing relates this definition with the existence of a replicating
portfolio for each bounded $\CF_T$--measurable random variable.

\begin{theorem}[FTAP II]
Suppose that assumption \ref{NFLVR} holds. Then the following are equivalent:
\begin{enumerate}
\item The market is complete (i.e. $\CM^e(S)=\{Q\}$).
\item For each $X\in L^\infty(\Omega,\CF_T,P)$ there exist a unique admissible portfolio $H\in\CH$ and a constant 
$x\in\RR$ such that
\be
X=x+(H\cdot S)_T.
\ee
\end{enumerate}
\label{FTAPII}
\end{theorem}

For complete markets, Merton's problem can be solved almost explicitly in terms of $\frac{dQ}{dP}$. The results of
the next two theorems, which are slightly stronger versions of \cite[theorem 2.0]{KraSch99} and 
\cite[theorem 2.1]{Scha01a}
(since we are assuming reasonable asymptotic elasticity for all our utility functions), are the analogues
of theorems \ref{case1} and \ref{case2} for complete markets. Notice that for either case 1 
or case 2 in assumption
\ref{utilass}, the value function for the dual problem is reduced to 
\be
v(y)=E\left[V\left(y\frac{dQ}{dP}\right)\right], \quad y>0
\label{vcomplete}
\ee
(for case 1 this was proved in \cite[lemma 4.3]{KraSch99}; for case 2 it is trivial, since 
$\CM^e(S)=\{Q\}$ implies that $\CM^a(S)=\{Q\}$ as well).

\begin{theorem} 
Suppose that $\CM^e(S)=\{Q\}$ and assumptions \ref{utilass} (case 1) and \ref{finite1} hold. Then, for any
$x\in\mbox{dom}(U)$, the problem 
\be
u(x)=\sup_{X_T\in C(x)} E[U(X_T)]
\ee
has a unique optimizer $\widehat{X}_T(x)\in C(x)$ given by
\be
\widehat{X}_T(x)=-V^\prime\left(y\frac{dQ}{dP}\right),
\ee
where $y$ is the solution to the equation
\be
E\left[-V^\prime\left(y\frac{dQ}{dP}\right)\frac{dQ}{dP}\right]=x.
\ee
\label{complete1}
\end{theorem}

\begin{theorem}
Suppose that $\CM^e(S)=\{Q\}$ and assumptions \ref{utilass} (case 2) and \ref{finite2} hold.
Then, for any $x\in\RR$, the problem  
\be
u(x)=\sup_{f\in C_U(x)} E[U(f)],
\ee
has a unique optimizer $\widehat{f}(x)\in C_U(x)$ given by
\be
\widehat{f}(x)=-V^\prime\left(y\frac{dQ}{dP}\right),
\ee
where $y$ is the solution to the equation
\be
E\left[-V^\prime\left(y\frac{dQ}{dP}\right)\frac{dQ}{dP}\right]=x.
\ee
Moreover, $\widehat{f}(x)$ equals the terminal value $\widehat{X}_T(x)$ of a uniformly integrable
$Q$--martingale of the form 
\[\widehat{X}_t(x)=x+(\widehat{H}(x)\cdot S)_t, \quad t\in[0,T],\]
for some $\widehat{H}(x)\in L(S)$.
\label{complete2}
\end{theorem}

There is no need to state versions of theorems \ref{hedge1} and \ref{hedge2}, since for complete markets the 
solution to the hedging problem \eqref{optimalhedge} for a bounded claim $B$ can be expressed in terms of
the solution of Merton's problem. Indeed, by theorem \ref{FTAPII}, there exists $(B_0,H^B)$ such that
\[B=B_0+(H^B\cdot S)_T\]
and this can now be used to write \eqref{optimalhedge} in the form of the Merton problem
\[\sup_{H\in\CH} E\left[U\left(x-B_0+((H-H^B)\cdot S)_T\right)\right].\]
Therefore, if $\widehat{H}^0(x-B_0)$ is the optimal portfolio for the Merton problem starting with
wealth $x-B_0$ obtained either from theorem \ref{complete1} or from theorem \ref{complete2}, then the optimal portfolio
for the hedging problem for the claim $B$ starting with wealth $x$ will be given by
\be
\widehat{H}(x)=\widehat{H}^0(x+B_0)+H^B.
\ee

\section{The dynamics of portfolio selection}

The theorems of the previous section give precise statements ensuring the existence and uniqueness of
solutions for both the optimal investment and optimal hedging problems for different types of utility functions. 
We have seen that under well defined conditions, there is a clear sense in which the optimal solution can always be 
approximated arbitrarily well by trading according to admissible portfolios. 
In what follows, to adopt a unified notation, we will write $H\in\CA$, which
loosely stands for ``allowed'' portfolios. In the back of our minds, however, we will keep the rigorous notion of what
it stands for:
admissible portfolios which, starting with initial capital $x\in\RR$, generate terminal wealths in $C(x)$ and $\CC(x)$ 
for theorems \ref{case1} and
\ref{hedge1}, respectively, or terminal wealths whose utilities arbitrarily approximate
the utility of the optimal solutions (in the $L^1$ sense) $\widehat{f}(x)$ for theorems \ref{case2} and 
\ref{hedge2}. We also use the notation $\CA_{(s,t]}$ for portfolio processes defined only on the time interval $(s,t]$,
as well as the shorthand for stochastic integration in this interval
\[(H\cdot S)_s^t:=\int_s^t H_udS_u , \quad 0\leq s\leq t\leq T.\]
Consistently with our previous section we have that $\CA=\CA_{(0,T]}$ and $(H\cdot S)_t=(H\cdot S)_0^t$, $t\in[0,T]$.

To understand better the optimal selection problem it is useful to formulate a dynamical version of it. Let us 
write $\widehat{H}^{(x,0)}$
for the optimal solution to the static primal problem
\be
u(x)= \sup_{H\in\CA} E[U(x+(H\cdot S)_T - B)],
\ee
obtained according to the theorems of the previous section, that is, starting at time $0$ with initial wealth $x$. 
For any intermediate
time $t\in[0,T]$ and $x\in\mbox{dom}(U)$, we can write
\begin{eqnarray}
u(x)&=& \sup_{H\in\CA} E[U(x+(H\cdot S)_T - B)] \nonumber \\
&=& \sup_{H\in\CA_{(0,t]}} E\left[\sup_{H\in\CA_{(t,T]}} E_t[U(x+(H\cdot S)_0^t+(H\cdot S)_t^T - B)]\right], 
\end{eqnarray}
which leads us to the study of the {\it conditional} problem
\begin{equation}
u_t({\sf w})=\sup_{H\in\CA_{(t,T]}} E_t[U({\sf w}+(H\cdot S)_t^T - B)],
\label{condprimal}
\end{equation}
where ${\sf w}\in\RR$ represents the wealth accumulated up to time $t$. If we trade according to $\widehat{H}^{(x,0)}$ up
to time $t$, that is, if ${\sf w}=x+(\widehat{H}^{(x,0)}\cdot S)_t$, then we must have
\begin{equation}
u_t({\sf w})=E_t[U({\sf w}+(\widehat{H}^{({\sf w},t)}\cdot S)_t^T - B)],
\end{equation}
for some portfolio $\widehat{H}^{({\sf w},t)}\in\CA_{(t,T]}$ (that is, starting at time $t$ with wealth ${\sf w}$)
which agrees with the restriction of $\widehat{H}^{(x,0)}$ on the 
interval $(t,T]$. In other words, the optimal portfolio $\widehat{H}^{(x,0)}$ is also conditionally optimal.
This is a special instance of the dynamic programming principle, which for this stochastic 
control problem has the form
\begin{equation}
u_s({\sf w})=\sup_{H\in\CA_{(s,t]}} E_s[u_t({\sf w}+(H\cdot S)_s^t)],
\label{dprog}
\end{equation}
for $0\leq s\leq t\leq T$.

\vspace{0.2in}
\noindent
{\bf The certainty equivalent value and the indifference price}
\vspace{0.2in}

There is a useful way to view the value function $u_t({\sf w})$. By the intermediate value theorem, $U^{-1}(u_t({\sf w}))$ exists
for each $({\sf w},t)$, $P$--almost surely. This defines, for each $({\sf w},t)$, the random variable
\be
B_t({\sf w})={\sf w}-U^{-1}(u_t({\sf w})),
\ee
which can be called the {\it certainty equivalent value} of the claim $B$ at time $t$. Since
\[U({\sf w}-B_t({\sf w}))=E_t[U({\sf w}+(\widehat{H}^{({\sf w},t)}\cdot S)_t^T - B)],\]
the certain utility achieved by investing the amount ${\sf w}-B_t({\sf w})$ in the risk free account equals the expected 
utility of the terminal
wealth ${\sf w}+(\widehat{H}^{({\sf w},t)}\cdot S)_t^T - B$ of the optimal hedging portfolio. 
For Merton's problem, where $B\equiv 0$,
the amount $-B^{0}_t({\sf w})$ indicates by how much the optimally invested portfolio outperforms the constant portfolio ${\sf w}$ 
over the period $(t,T]$. By putting $s=0$ in \eqref{dprog} we obtain
\begin{eqnarray}
u(x)&=&\sup_{H\in\CA_{(0,t]}} E[u_t(x+(H\cdot S)_0^t] \nonumber \\
&=&  \sup_{H\in\CA_{(0,t]}} E[U[x+(H\cdot S)_0^t-B_t(x+(H\cdot S)_0^t)].
\label{effective}
\end{eqnarray}
Therefore, $B_t({\sf w})$ represents a wealth dependent effective value of the claim $B$ at time $t$.

Following \cite{HodNeu89} (according to \cite{Bech01}), a clear interpretation of the certainty equivalent values
can be given by considering an investor who, holding wealth ${\sf w}$ at time
$t$, must decide the minimum amount $\pi$ to charge when selling a claim $B$. If he sells the claim for 
$\pi$ and hedges optimally against the claim by holding the portfolio $\widehat{H}^{({\sf w}+\pi,t)}$, he will 
achieve an expected utility
\[E_t[U({\sf w}+\pi+(\widehat{H}^{({\sf w}+\pi,t)}\cdot S)_t^T - B)]=U({\sf w}+\pi-B_t({\sf w}+\pi))\]
If, however, he does not sell the claim and invests optimally for Merton's problem, he achieves 
\[E_t[U({\sf w}+(\widehat{H}^{({\sf w},t)}(0)\cdot S)_t^T)]=U({\sf w}-B^0_t({\sf w})).\]
The {\it indifference price} of the claim $B$ at time $t$ for wealth ${\sf w}$ is the value 
for $\pi=\pi_t^{B}({\sf w})$ which makes these equal, that is, it is the solution of
\be
\pi^B_t({\sf w})=B_t({\sf w}+\pi^B_t({\sf w}))-B^0_t({\sf w}).
\label{indiffsell}
\ee
Since we have defined these concepts from the point of view of an agent faced with a liability $B$ at time
$t$, this indifference price corresponds to a ``seller's price''. To obtain the correct notion of a ``buyer's price'',
we just need to consider the reverse claim $-B$, which then produces a terminal wealth with expected utility equaling that
of ${\sf w}-\pi-B_t({\sf w}-\pi)$ when bought by $\pi$. The indifference price is now the value of $\pi$ that makes this equal to
the amount whose certain utility equals the expected utility for Merton's problem starting with wealth ${\sf w}$ at time $t$,
which by definition is ${\sf w}-B^0_t({\sf w})$. In other words, it is the solution of
\be
\pi=B^0_t({\sf w})-B_t({\sf w}-\pi),
\label{indiffbuy}
\ee
which therefore equals $-\pi^B_t({\sf w})$ as defined in \eqref{indiffsell}.

In a complete market, the indifference price equals the risk-neutral price, that is, if the bounded claim
$B$ is written in terms of its unique replicating admissible portfolio $H^B$ as $B=B_0+(H^B\cdot S)_T$, then
\be
\pi^B_t=B_0+(H^B\cdot S)_0^t=E_{t,Q}[B],
\ee
where $Q$ is the unique equivalent martingale measure. The first equality above remains true in incomplete markets
if the claim $B$ happens to satisfy $B=B_0+(H^B\cdot S)_T$ for some admissible portfolio $H^B$. 

\vspace{0.2in}
\noindent
{\bf The Davis price}
\vspace{0.2in}

Let us assume for a moment  that the solutions of the dual problems in theorems \ref{case1} and \ref{case2} are
equivalent martingale measures (in case 2 we have seen that this indeed the case for the exponential utility under
the finite entropy condition; for counterexamples where in case 1 the solution fails to be a martingale, see 
\cite{KraSch99}). If, for each $\varepsilon\geq 0$, we let $B_t^\varepsilon ({\sf w})$ denote the 
certainty equivalent value of $\varepsilon B$,
then the Davis price of $B$ is defined to be \cite{Davi97}
\be
\pi^{\scriptscriptstyle Davis}_t({\sf w})=\left.\frac{dB_t^\varepsilon ({\sf w})}{d\varepsilon}\right|_{\varepsilon=0}.
\label{davis}
\ee

By differentiating the identity
\[U({\sf w}-B_t^\varepsilon ({\sf w}))=E_t[U({\sf w}+(\widehat{H}^{(\varepsilon,{\sf w},t)}\cdot S)_t^T 
- \varepsilon B)]\]
at $\varepsilon=0$ and noting that, by optimality, 
\[\left.\frac{d\widehat{H}^{(\varepsilon,{\sf w},t)}}{d\varepsilon}\right|_{\varepsilon=0}=0,\]
we see that
\[\pi^{\scriptscriptstyle Davis}_t({\sf w})=\frac{E_t[U^\prime({\sf w}+(\widehat{H}^{(0,{\sf w},t)}\cdot S)_t^T)B]}
{U^\prime({\sf w}-B^0_t({\sf w}))}.\]
But from the theory of the Merton problem, a dynamical version of either \eqref{rel1} or \eqref{rel2} gives
\[U^\prime({\sf w}+(\widehat{H}^{(0,{\sf w},t)}\cdot S)_t^T)=U^\prime({\sf w}-B^0_t({\sf w}))\frac{d\widehat{Q}_t(y)}{dP}\]
for $y=u_t^\prime({\sf w})$, where $\widehat{Q}_t(y)$ stands for the optimal solution to the conditional dual problem.
Thus the Davis price of $B$ is given by the expectation pricing 
\be
\pi^{\scriptscriptstyle Davis}_t({\sf w})=E_{t,\widehat{Q}_t(y)}[B].
\ee

We remark that the indifference price, being intrinsically nonlinear, does not in general satisfy useful criteria such 
as put-call parity. The Davis price, on the other hand, does.

\vspace{0.2in}
\noindent
{\bf Exponential utility}
\vspace{0.2in}

An important simplification occurs if we specialize to the exponential utility $U(x)=-\frac{e^{-\gamma x}}{\gamma}$, 
$\gamma>0$. A look
at \eqref{condprimal} shows that $u_t({\sf w})$ factorizes as
\be
u_t({\sf w})=-\frac{e^{-\gamma {\sf w}}}{\gamma}\inf_{H\in\CA_{(t,T]}} E_t\left[e^{-\gamma(H\cdot S)_t^T+\gamma B}\right]
=:-\frac{e^{-\gamma {\sf w}}}{\gamma}v_t.
\label{voptimal}
\ee
Here we see that $v_t$ is a time dependent but wealth independent $\CF_t$--measurable random variable. We also see
that the certainty equivalent value
\be
B_t=-\frac{1}{\gamma}\log v_t,
\ee
the optimal portfolio $\widehat{H}^{(t)}$ and the indifference price $\pi_t$ are all wealth independent processes.

\vspace{0.2in}
\noindent
{\bf Example (geometric Brownian motion):} Consider now a market of $d$ stocks whose prices, 
discounted by the constant interest rate $r$, satisfy
\be
\frac{dS^i_t}{S^i_t}=(\mu^i-r)dt+\sum_{\alpha=1}^d \sigma^{i\alpha} dW^\alpha,
\label{gbm}
\ee
where $\mu^i\in\RR$ and the invertible $d\times d$ matrix $\sigma^{i\alpha}$ are constants and $(W^\alpha)$ is
a $d$--dimensional $P$--Brownian motion. This market is complete and, as is well known, the unique
equivalent martingale measure $Q$ has Radon-Nikodym derivative
\be
\frac{dQ}{dP}=\exp\left[-\int_0^T\left(\sum_{\alpha} \lambda^\alpha dW^\alpha + \frac{1}{2}\|\lambda\|^2dt\right)\right],
\label{riskneutral}
\ee
with constant market price of risk $\lambda^\alpha = \sum_i (\sigma^{-1})^{\alpha i}(\mu^i-r)$.

For the exponential utility function with initial wealth $x$, the optimal discounted terminal 
wealth $\widehat{X}_T$ is given by
\be
e^{-\gamma\widehat{X}_T}=y\frac{dQ}{dP},
\ee
for a constant $y$ to be determined. From this and \eqref{riskneutral}, one finds
\[\widehat{X}_T=-\frac{1}{\gamma}\left(\log y + \frac{1}{2}\int_0^T\|\lambda\|^2dt\right)+
\frac{1}{\gamma}\int_0^T\sum_{ij}(\mu^i-r)\left((\sigma\sigma^T)^{-1}\right)^{ij}\frac{dS^j_t}{S^j_t},\]
so the optimal portfolio for Merton's problem is
\be
(\widehat{H}^0)^j_t=\frac{\sum_i \left((\sigma\sigma^T)^{-1}\right)^{ij}(\mu^i-r)}{\gamma S^j_t}
\label{optimalgbm}
\ee
and $y$ is the solution to the equation
\be
x=-\frac{1}{\gamma}\left(\log y +\frac{1}{2}\|\lambda\|^2T\right).
\label{ygbm}
\ee
The certainty equivalent value for Merton's problem in this market turns out to be
\be
B^0_t=\frac{1}{2}\|\lambda\|^2(t-T).
\ee

Since the market is complete, the indifference price of any bounded claim $B$ equals its risk-neutral
price (its Black-Scholes price), so the certainty equivalent value is given by
\be
B_t=\frac{1}{2}\|\lambda\|^2(t-T) + E_{t,Q}[B].
\ee
Finally, the optimal hedging portfolio for the claim $B$ is
\be
\widehat{H}_t=\widehat{H}^0_t+H^B_t,
\label{hedgegbm}
\ee
where $H^B_t$ is the replicating portfolio for $B$.

\section{Discrete time hedging}

We now restrict to discrete time hedging, where the portfolio
processes have the form
\begin{equation}
H_t=\sum^{K}_{k=1}\ H_k\ {\bf 1}_{(t_{k-1},t_{k}]}(t)
\label{discretehedge}
\end{equation}
where each $H_k$ is an $\RR^d$--valued $\CF_{k-1}$ random variable.
We take the discrete time partition of the interval $[0,T]$ to be of the form
\[t_0=0< t_1=\frac{T}{K} < \ldots < t_k=\frac{kT}{K} \ldots <t_K=T\]
and use the notation $S_j:= S_{t_j}$ for discrete time stochastic processes. The
discounted wealth process will be
$X_{j}=x+(H\cdot S)_j$, with the notation  
$(H\cdot S)_j:= \sum_{k=1}^{j} H_k\Delta S_k$, $ (H\cdot
S)^j_k:=(H\cdot S)_j-(H\cdot S)_k$ and
$\Delta S_k:=S_{{k}}-S_{k-1}$.

Now the dynamic programming problem (\ref{dprog}) falls into $K$
subproblems
\begin{equation}
u_{k-1}(x)=\sup_{H\in \CA_{(t_{k-1},t_k]}}E_{k-1}[u_k(x+H_k\Delta
S_k)],\qquad k=K,K-1,\dots,1
\label{uiteration}
\end{equation}
subject to the terminal condition $u_K(x)=U(x)$. Then, for each $x$ this
defines a process
$u_k(x)$. Similarly, the certainty equivalent value process $B_k(x)$ is
defined iteratively by
\begin{equation}
U(x-B_{k-1}(x))=\sup_{H\in \CA_{(t_{k-1},t_k]}}E_{k-1}[U(x+H_k\Delta
S_k -B_k(x+H_k\Delta S_k)]
\label{Biteration}
\end{equation}
with $B_K(x)$ taken equal to the terminal claim $B$. In both formulations
of the problem, let
$\widehat H_k(x)$ denote the minimizer, which of course is an $\CF_{k-1}$ random variable.

In what follows, we will be largely concerned with markets and claims which satisfy
the Markovian conditions:

\begin{assumption}
\label{Markovmarket} The market is Markovian and its state variables\\
$Z=(S^1,\ldots, S^d,Y^1,\ldots, Y^{n-d})$ lie in a finite dimensional
state space $\CS\in \RR^n$.
\end{assumption}

\begin{assumption} The contingent claim is taken to be of the form $B_T=\Phi(Z_T)$
for a bounded Borel function $\Phi:\CS\to \RR$.
\label{Markovclaim}
\end{assumption}

In these assumptions, we interpret $S$ as discounted asset prices as before
and the additional variables
$Y$ as values of nontraded quantities such as stochastic volatilities
which may or may not be observed directly. In this Markovian setting, the solution of (\ref{uiteration}) and 
the optimal allocation have the form
\begin{eqnarray}
u_k(x)&=&g_k(x,Z_k)\\
\widehat H_{k+1}(x)&=& h_{k+1}(x,Z_{k})\label{hstrat}
\end{eqnarray}
for (deterministic) Borel functions
$\{g_k,h_{k+1}\}_{k=0}^{K-1}$ mapping $\mbox{dom}(U)\times \CS$ to $\RR$ and
$\RR^d$ respectively. Similarly, the solution of (\ref{Biteration}) has the optimal allocation
$\widehat H_{k+1}$ as above and $B_k$ in the form
\begin{equation}
B_k(x)=b_k(x,Z_k)
\label{Bsolution}
\end{equation}
for  Borel functions
$\{b_k\}_{k=0}^{K-1}$ mapping $\mbox{dom}(U)\times \CS$ to $\RR$.

As indicated in the previous section, matters simplify in the special case
of the exponential utility function
$U(x)=-\frac{e^{-\gamma x}}{\gamma}$, $\gamma>0$. One finds that the dynamic 
program
can be written in the wealth independent form $u_k(x)=U(x)v_k$, $\widehat H_k(x)=\widehat H_k$,
and $B_k(x)=B_k$ where  the random variables $v$, $\widehat H$ and $B$
have the form:
\begin{eqnarray}
v_k&=&g_k(Z_k)\\
\widehat H_{k+1}&=& h_{k+1}(Z_{k})\label{exphstrat} \\
B_k &=& b_k(Z_k)
\end{eqnarray}
for deterministic functions $g_k,h_{k+1}$, and $b_k$ on the state space $\CS$.
The iteration equations are simply
\begin{eqnarray}
g_k(Z)&=& \inf_{h\in \RR^d}E_k[\exp(-\gamma h\cdot\Delta S_{k+1})
g_{k+1}(Z_{k+1}) | Z_k=Z]  \\
\exp(\gamma b_k(Z)) &=&\inf_{h\in \RR^d}E_k[\exp(-\gamma(h\cdot\Delta
S_{k+1} - b_{k+1}(Z_{k+1})) | Z_k=Z]
\end{eqnarray}
and the optimal $h$ defines the function $h_{k+1}(Z)$.

\section{The exponential utility allocation algorithm}

In this section we introduce a Monte Carlo method for learning the
optimal trading strategy (\ref{hstrat}) for the discrete time Markovian
problems discussed in the previous section. We want an algorithm which
will generate an approximate trading rule, based on a data set
$\{Z^i_k\}_{i=1,\dots,N;k=0,\dots,K}$ where $Z^i_k\in \RR^n$ denotes
the ``state '' of the $i$th sample path at time $t_k=kT/K$.

Consider the discrete time problem \eqref{uiteration} for a general utility
function $U$ satisfying assumption \ref{utilass}. The optimal portfolio $\widehat{H}^{i}_{k+1}\in \RR^d$
should be selected as $h_{k+1}(X^{i}_k,Z^{i}_k)$ where $X^i_k$ is the
wealth held at the point $(i,k)$. We can see a basic difficulty with a Monte
Carlo approach: to ``learn'' the function $h_{k+1}$ from the data
$\{Z\}$ will require being able to fill in the optimal wealth from
time $t_0$ to $t_k$. Then, conditionally upon knowing the wealth $X^i_k$
at time
$t_k$, finding the function $h_{k+1}$ requires dynamic programming
backwards from time
$t_K=T$ to $t_k$. In other words, a Monte Carlo learning algorithm for
$\widehat{H}^{i}_k$ based on general utility will require both forward and backward
dynamic programming. We have no effective method to suggest for general
utility.

By contrast, in the special case of exponential utility, the theoretical
optimal rule $\widehat{H}^{i}_{k+1}=h_j(Z^{i}_k)$ depends only on the directly observed
data $\{Z^{i}_k\}$ and is independent of the wealth $X^{i}_k$. For this
reason our algorithm works only for exponential utility, and we take
\[ U(x)=-e^{-x},\]
for simplicity.

\subsection{The algorithm}

\begin{enumerate}
\item Step  $k=K$: The final optimal allocation $\widehat H_K$ is defined to
be the $\RR^d$--valued $\CF_{K-1}$ random variable which solves
\[ \min_{H\in \CA_{(K-1,K]}} E_{K-1}[\exp(-H\cdot\Delta S_{K}
+B)],\qquad \mbox{(a.s.)}
\]
which is easily seen to be equivalent to the optimizer of
\begin{equation}
\min_{H\in \CA_{(K-1,K]}} E[\exp(-H\cdot\Delta S_{K}
+B)]
\end{equation}
Since the solution is known to be given by $\widehat H_{K}=h_K(Z_{K-1})$ for
some deterministic function $h_K\in\CB(\CS)$ (the set of Borel functions on $\CS$), we write this as
\begin{equation}
\min_{h\in\CC(\CS)} \Psi_{K}(h)
\end{equation}
where $\Psi_{K}(h):=E[\exp(-h(Z_{K-1})\cdot\Delta S_{K}
+B)]$. On a finite set of data, we can pick an $R$--dimensional subspace
$\CR(\CS)\subset \CB(\CS)$ of functions on $\CS$ and attempt to  ``learn'' a suboptimal solution
\[h^\CR_K =\argmin_{h\in\CR(\CS)}\Psi_{K}(h)\]
By the central limit theorem, the expectation $\Psi_{K-1}(h)$ for $h$ in
a neighbourhood of $h^\CR_K$, and hence the solution $h^\CR_K$ itself, 
can
be approximated by the finite sample estimate
\begin{equation}
\widetilde{\Psi}_{K}(h)=\frac1{N}\sum_{i=1}^N\
\exp\left(-h(Z^i_{K-1})\cdot\Delta S^i_K +\Phi(Z^i_K)\right)
\end{equation}
This leads to the estimator $\widetilde{h}^\CR_{K}$ based on $\{Z^{i}_k\}$
and the choice of subspace $\CR$ defined by
\begin{equation}
\label{finiteoptimization}
\widetilde{h}^\CR_{K} = \argmin_{h\in\CR(\CS)}\widetilde{\Psi}_{K}(h)
\end{equation}

\item Inductive step for $k=K-1,\dots,2$:  The estimate $\widetilde{h}^\CR_k$ of the optimal 
rule $\widehat h_k$, for $2\le k <K-1$ is
determined inductively given the estimates $\widetilde{h}^\CR_{k+1},\dots,\widetilde{h}^\CR_{K}$. It is defined to be
\begin{equation}
\widetilde{h}^\CR_k = \argmin_{h\in\CR(\CS)}\widetilde{\Psi}_{k}(h;\widetilde{h}^\CR_{k+1},\dots,
\widetilde{h}^\CR_{K} )
\end{equation}
where
\begin{eqnarray}
\widetilde{\Psi}_{k}(h;\widetilde{h}^\CR_{k+1},\dots,\widetilde{h}^\CR_{K} ) &=&\nonumber\\
&&\hspace{-2.1in}
\frac1{N}\sum_{i=1}^N\
\exp\left(-h(Z^i_{k})\cdot\Delta S^i_{k+1}-\sum_{j=k+1}^{K} \widetilde{h}^\CR_j(Z^i_j)
\cdot\Delta S^i_{j} +\Phi( Z^i_K)\right)
\end{eqnarray}

\item Final step $k=1$:  This step is degenerate since the initial
values $Z_0$ are constant over the sample. Therefore we determine the
optimal constant vector $\widetilde{h}_1\in\RR^d$ by solving
\begin{equation}
\widetilde{h}_1 = \argmin_{h\in\RR^d}\widetilde{\Psi}_{1}(h;\widetilde{h}^\CR_{2},\dots,
\widetilde{h}^\CR_{K} )
\end{equation}
\end{enumerate}

To summarize, the algorithm above learns a collection of functions of the form
$(\widetilde{h}_1,\widetilde{h}^\CR_2,\dots,\widetilde{h}^\CR_{K})\in \RR^d \times
\CR(\CS)^{K-1}$ from the Monte Carlo simulation. This collection defines a
suboptimal allocation strategy for the exponential hedging problem. Finally, the optimal value
$\widetilde{\Psi}_{1}(\widetilde{h}_1;\widetilde{h}^\CR_{2},\dots,\widetilde{h}^\CR_{K})$ is an estimate of the quantity
$e^{B_0}$, where $B_0$ is the certainty equivalent value of the claim $B$ at time $t=0$.

\subsection{Discussion of errors}

It is important to identify two distinct systematic sources of error in
the algorithm. The first, which we call {\it approximation one}, is in
focusing on suboptimal solutions $h^\CR_k$ which lie in a specified
subspace
$\CR(\CS)$ of the full space $\CB(\CS)$. From a pragmatic perspective, we need to
select a set of $R$ basis functions $f_1,\dots,f_R$ for $\CR(\CS)$ which
does a good job of representing the true optimal function over the values
of state space covered by the Monte Carlo simulation. Naively, one might
expect to need to choose $R$ exponentially related to the dimension of
$\CS$; experience seems to indicate far fewer functions are needed for
higher dimension problems. For a discussion of this type of question in
the context of the Longstaff-Schwartz (LS) method for American options, see
\cite{LonSch01} and \cite{CleLamPro02}. Observe that the requirements of
our algorithm are much more stringent than for the American
option problem, since the strategy to be learned is not simply ``to exercise or
not to exercise'', but must select a high dimensional vector at
each point $(i,k)$ in the simulation. Having said this, we take the point of view that the careful
selection of a subspace
$\CR(\CS)$ might lead to good performance of the algorithm. Furthermore, our
experiments show that the sensitivity  to changes in $\CR(\CS)$ of quantities
such as indifference prices are much less than that of quantities such as
hedge allocations.

The second source of error, {\it approximation two}, is the finite $N$
approximation. We can in principle estimate this error in terms of the basic
model parameters; the following is a heuristic argument to give a
flavour of the problem for the $k$th step, $k\le K$.  By the central
limit theorem, for a given confidence level
$1-\alpha,
\alpha\ll 1$, there exist constants
$C_1,C_2$ so that
\[ \|\Psi(h)-\widetilde\Psi(h)\|\le
\frac{C_1}{\sqrt{N}},
\quad\|\nabla\Psi(h)-\nabla{\widetilde\Psi}(h)\|\le\frac{C_2}{\sqrt{N}}
\]
with probability $1-\alpha$, for $h$ in a convex neighbourhood of the
true critical point $h^{\CR}$, defined by $\nabla\Psi(h^\CR)=0$. We suppose that the estimated 
critical point $\widetilde h^\CR$, defined by $\nabla{\widetilde\Psi}(\widetilde h^\CR)=0$, lies in this
neighbourhood, and furthermore the operator inequalities 
$0<C_3\leq\nabla^2\Psi\leq C_4$
hold on the same neighbourhood. Then one immediately 
derives
the inequalities
\begin{eqnarray}
\|h^\CR-\widetilde h^\CR\|&\le& \frac{C_2}{C_3\sqrt{N}} \label{error1}\\
|\Psi(h)-\widetilde\Psi(\widetilde h^\CR)|&\le&
\frac{C_1}{\sqrt{N}}+\frac{C_2^2C_4}{2C_3^2N}
\label{error2}
\end{eqnarray}
which show convergence of $\widetilde h^\CR$ to $h^\CR$ for  $N\to \infty$.

The above discussion addresses the errors made at the $k$th time step of
the algorithm. Further study is needed to understand how errors
accumulate as $k$ is iterated. The answer to this question will give
guidance on how to distribute computational effort over the different
time steps, and can be expected to parallel the same question as it
arises for the LS algorithm.

\section{Numerical implementation}

We have tested the algorithm in the
simple problem of  investment with exponential utility $U(x)=-e^{-x}$ 
in a
stock which behaves as a geometric Brownian motion, with and without the
purchase of a single at--the--money European put  option. We have seen
there is an exact solution to this problem which can be compared in
detail to the solution generated by the Monte Carlo algorithm.

We consider the model of \eqref{gbm} with $d=1$ and parameters
$S_0=1$, $\mu=0.1$, $\sigma=0.2$ and $r=0.0$ over the period of one year $T=1$. We
apply the allocation algorithm to two scenarios involving portfolio
selection at discrete time intervals of $1/50$ (i.e. weekly): i) the
Merton investment problem; and ii)  the hedging problem for the buyer of a single
written at--the--money European put. In each case we apply the method for
simulations of length
$N=1000,10000,100000$. Then, for comparison to theory, we use the same
Monte Carlo simulations, but rehedged
weekly according to the theoretical formula \eqref{hedgegbm}, with $-H^B$ equal to the 
Black-Scholes delta of the option.

Our results are displayed in figures 1 to 5. Figures 1,2,3,4 show the
profit/loss distributions at time $T$ for the learned Merton, learned put
option, true Merton and true put option cases respectively. They show the
empirical distributions for $N=1000,10000,100000$. Figure 5 shows the
values of the hedge ratio along a single sample path calculated according
to both the strategy learned with $N=100000$ and the true strategy.

For comparison of their performances, we tabulate below the mean and the standard deviation of these
distributions in each of the four cases,
as well as the final expected exponential utility with parameters $\gamma=1/4$ $(U_1)$, $\gamma=1$ $(U_2)$ and
$\gamma=4$ $(U_3)$, corresponding to an increasing order of risk-aversion. As measures of the risk
associated with each case, we also tabulate their value-at-risk and conditional value-at-risk 
for $90\%$ ($\mbox{VaR}_{90}$ and $\mbox{CVar}_{90}$) and $99\%$ 
($\mbox{VaR}_{99}$ and $\mbox{CVar}_{99}$) confidence levels. 

From the $U_2$ values on the table, one can derive the learned estimates of the indifference price
$0.0767,0.0790$ and $0.0792$ for the cases $N=1000,10000$ and $100000$. Using the true strategy leads to the values
$0.0798,0.0796,0.0795$, respectively. The theoretical Black-Scholes price is $0.0797$.  

\begin{table}[h]
\hspace{-0.6in}\begin{tabular}{|c|c|c|c|c|c|c|c|c|c|}
\hline
Case& Mean& St. Dev. &$U_1$ & $U_2$ & $U_3$
&$\mbox{VaR}_{99}$&$\mbox{VaR}_{90}$&$\mbox{CVaR}_{99}$&$\mbox{CVaR}_{90}$\\
\hline\hline
1a&    0.3572   & 0.5778  & -0.9241   & -0.8262   & -3.3023 &  -1.0829 &
-0.3659 &  -1.2053  & -0.6495

\\\hline
1b&    0.2768 &   0.5136   & -0.9409  &  -0.8674  &  -3.5224 &  -1.0159
&-0.3743
  &   -1.2144 &  -0.6543

\\\hline
1c &   0.2528 &   0.5013  &  -0.9462 &   -0.8810  &  -2.7684 &  -0.9209
&-0.3913 &   -1.0913 &  -0.6318

\\\hline
\hline
2a&    0.4349   & 0.5797  &  -0.9064   & -0.7652  &  -2.4748  & -0.9828 &
-0.3012 &  -1.1348  & -0.5720

\\\hline
2b&    0.3562   & 0.5142   & -0.9224 &   -0.8015  &  -2.9355  & -0.9518  &
-0.2732
  &  -1.1756   &-0.5626

\\\hline

2c&    0.3325 &   0.5020 &   -0.9275 &   -0.8139 &   -2.4439 &  -0.8532 &
-0.2859&   -1.0633 &  -0.5387\\\hline
\hline
3a&    0.2307   & 0.4898 &   -0.9511 &   -0.8956  &  -2.5283  & -0.9723 &
-0.3961 &    -1.0318  & -0.6430

\\\hline
  3b&  0.2524  &  0.4945 &   -0.9460  &  -0.8773   & -2.4184 &  -0.8852  &
-0.3778 &   -1.0429 &  -0.6096

\\\hline

3c&    0.2506  &  0.4995  &  -0.9466  &  -0.8816  &  -2.6461&   -0.9081&
-0.3896 &   -1.0733  & -0.6254

\\\hline
\hline
4a&    0.3108  &  0.4904  &  -0.9322   & -0.8269  &  -1.8302 &  -0.8878
&-0.3135 &   -0.9492 &  -0.5628

\\\hline
4b&    0.3322 &   0.4954   & -0.9274 &   -0.8102  &  -1.7521 &  -0.8054 &
-0.2979 &   -0.9616  & -0.5289

\\\hline
4c&    0.3304  &  0.5005 &   -0.9279  &  -0.8142 &   -1.9178 &  -0.8274
&  -0.3098
  &    -0.9921  & -0.5449

\\\hline
\end{tabular}
\caption{Mean, standard deviation, final expected utilities and risk measures for the
profit/loss distribution of the learned Merton, learned put option, true Merton and true put option
portfolios with (a) 1000, (b) 10000 and (c) 100000 Monte Carlo simulations of stock prices following a
geometric Brownian motion.}
\label{risktable}
\end{table}

\section{Discussion}

This paper seeks to bridge the gap between the theory of exponential 
hedging in incomplete markets and the numerical implementation of that 
theory. Utility based hedging introduces several key concepts, notably 
certainty equivalent values and indifference prices which have no 
counterpart in complete markets. Therefore we have little experience or 
intuition on which to base our understanding of optimal trading in these 
markets. The simple and flexible Monte Carlo algorithm we introduce in 
this paper provides a test bed for realizing the theory of exponential 
hedging in essentially any market model. For example, problems involving 
American style early-exercise options can in principle be easily 
included in our framework by following the Longstaff-Schwartz Monte 
Carlo method \cite{LonSch01}. Using our method for a variety 
of problems should help one gain intuition and understanding of how 
exponential hedging works in practice and how it compares with other 
hedging approaches.

Our preliminary study of the geometric Brownian model shows not 
unexpectedly that the method performs better for pricing than hedging. 
Interestingly the indifference price, perhaps the key theoretical 
concept, appears to be better approximated than the two certainty 
equivalent values which define it. On the other hand, as we see from the 
sample path shown in figure \ref{hedgeplot}, the actual hedging 
strategy learned by the algorithm deviates a lot from the theoretical  
strategy along individual stock trajectories, and cannot be seen as 
reliable.

Predictably, the basic method we use shows some distinct shortcomings 
which prevent it from being taken as a {\it de jure} guide to real trading. 
Approximation one arises by restricting possible hedge strategies to a 
low dimensional subspace. It is clear that such a restriction will often 
lead to unsuitable strategies. However, we feel that approximation two, 
the finite sample size error, will likely be even more problematic for 
practical realizations. A brief study of the size of  the constants 
which enter the estimates (\ref{error1}) and (\ref{error2}) suggests that reliable  
learned strategies will demand a very large value of $N$. In 
our simulations, $N=100000$ gave reliable prices, but not hedging 
strategies.  A third difficulty we noticed arising in our method is that 
learned strategies fluctuate far too much in time. Some simple smoothing 
procedure in time might lead to a marked improvement in hedging.

To conclude this discussion, it is worthwhile to revisit the way in 
which our method of dynamic programming (finding $\widehat H$ by induction over $K$ steps 
backwards in time) leads to computational efficiency compared to a more 
direct approach which seeks to compute the optimal hedging strategy 
$\widehat H$  simultaneously at all times. Fixing as before an $R$--dimensional subspace $\CR(\CS)$  
for the form of the hedging strategy at each time, direct optimization 
of a single convex function of $K\times R$ variables costs 
$\CO(NR^2K^2)$ flops. By dynamic programming this is reduced to $K$  
sequential optimizations of functions of $R$ variables which will take 
$\CO(NR^2K)$ flops. Accuracy is preserved by dynamic programming 
because the $KR\times KR$ Hessian matrix of the global optimization is 
approximately block diagonal over the individual time steps.

Putting aside the obvious drawbacks of the algorithm, we can see that 
our very simple and direct method will shed light on most conceptual 
difficulties arising in exponential hedging in incomplete markets. It 
implements the spirit of dynamic programming and prices claims quite 
reliably, even if it cannot easily produce accurate estimates of hedging 
strategies. On these merits alone, we think our algorithm deserves much 
further study and refinement.

\newpage
\begin{figure}[ht]
\epsfxsize=360pt
\epsfysize=500pt
\centerline{\epsffile{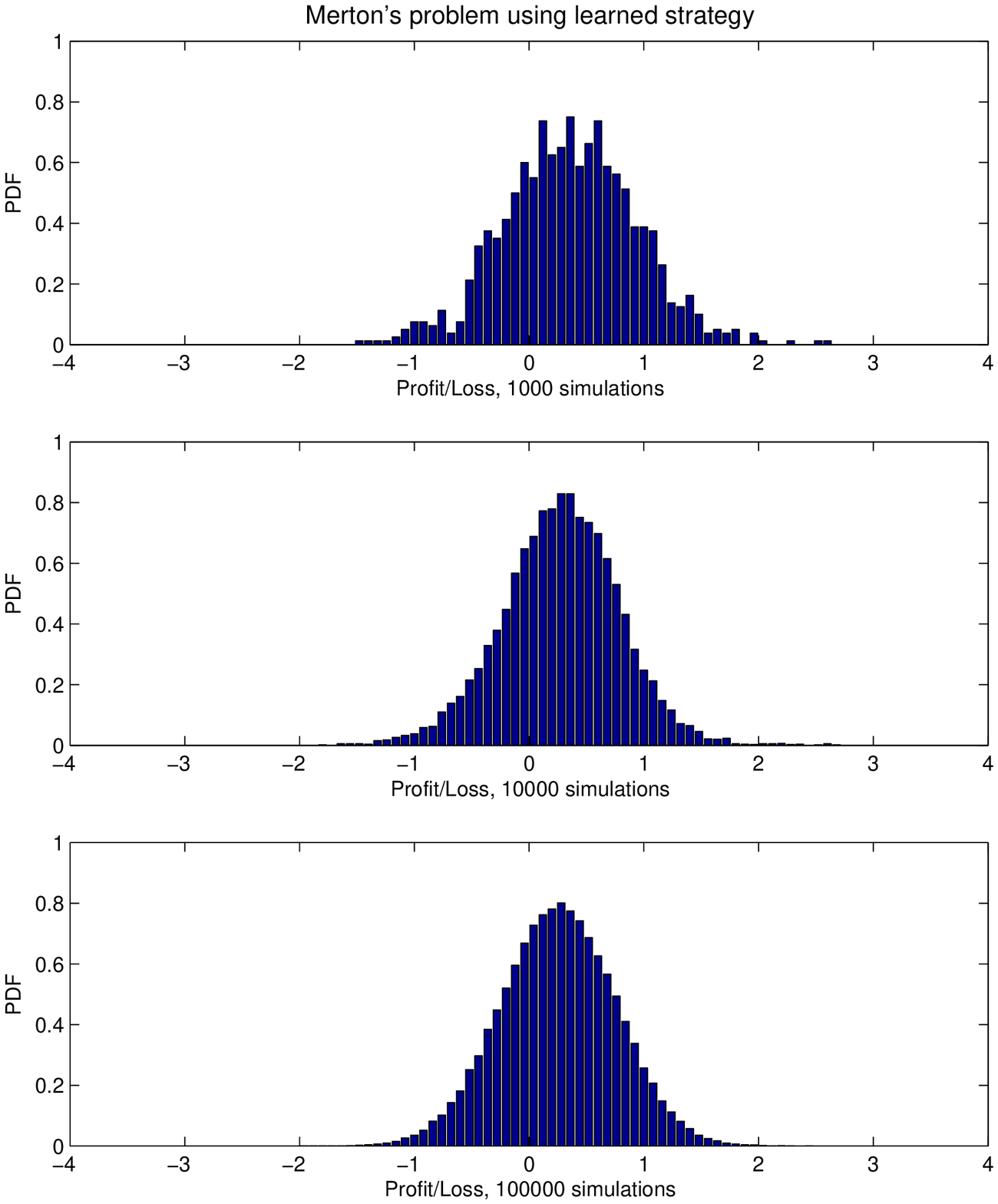}}
\caption{The profit/loss distribution of the learned investment portfolio,
obtained from the exponential utility allocation algorithm as an approximated solution to Merton's problem,
evaluated on simulated stock prices following a 
geometric Brownian motion.}
\label{mertonfig}
\end{figure}

\newpage
\begin{figure}[ht]
\epsfxsize=360pt
\epsfysize=500pt
\centerline{\epsffile{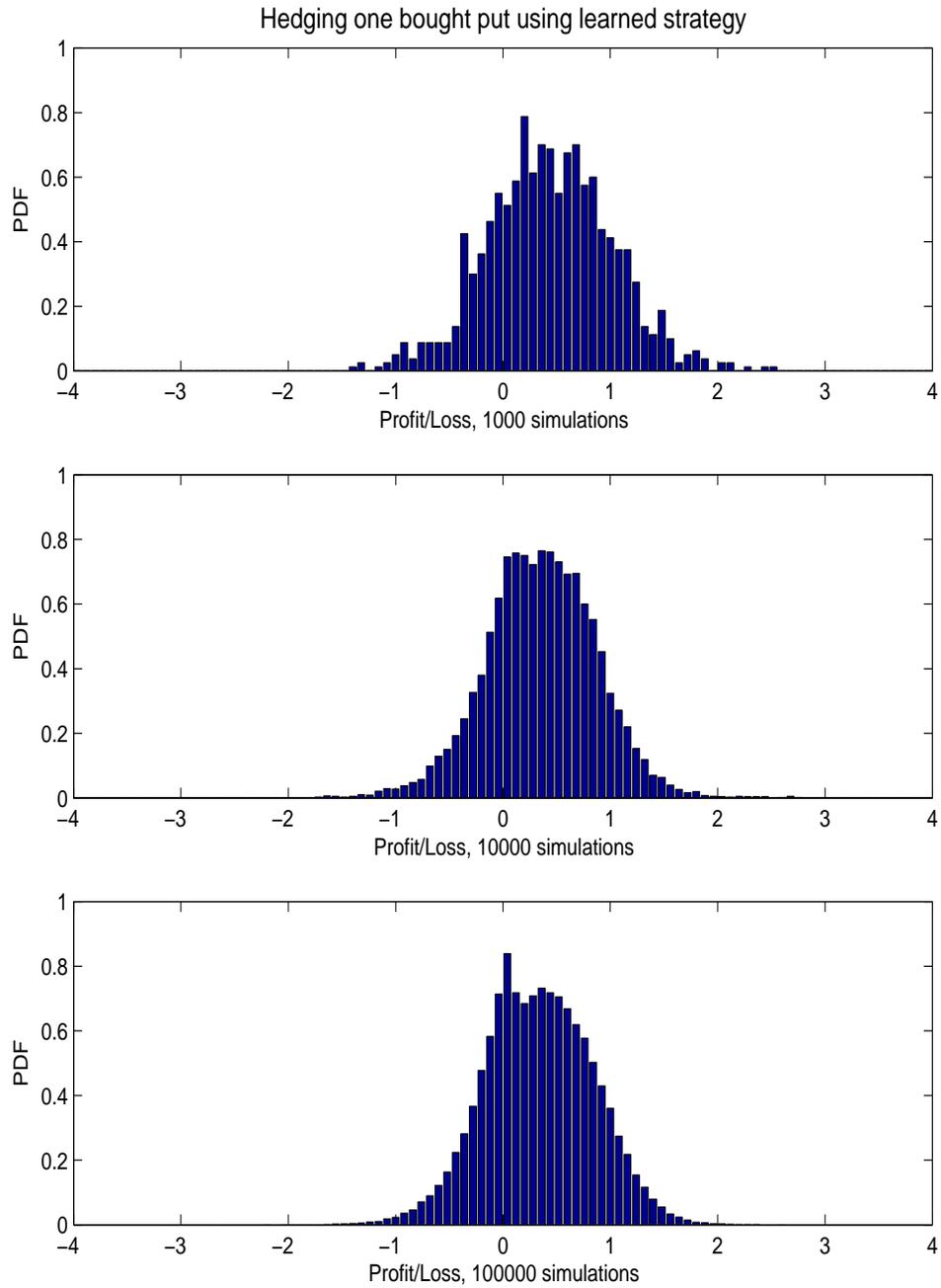}}
\caption{The profit/loss distribution of the learned hedging portfolio for the buyer of one put option,
obtained from the exponential utility allocation algorithm on simulated stock prices following a 
geometric Brownian motion.}
\label{oneputfig}
\end{figure}

\newpage
\epsfxsize=360pt
\epsfysize=500pt
\begin{figure}[ht]
\centerline{\epsffile{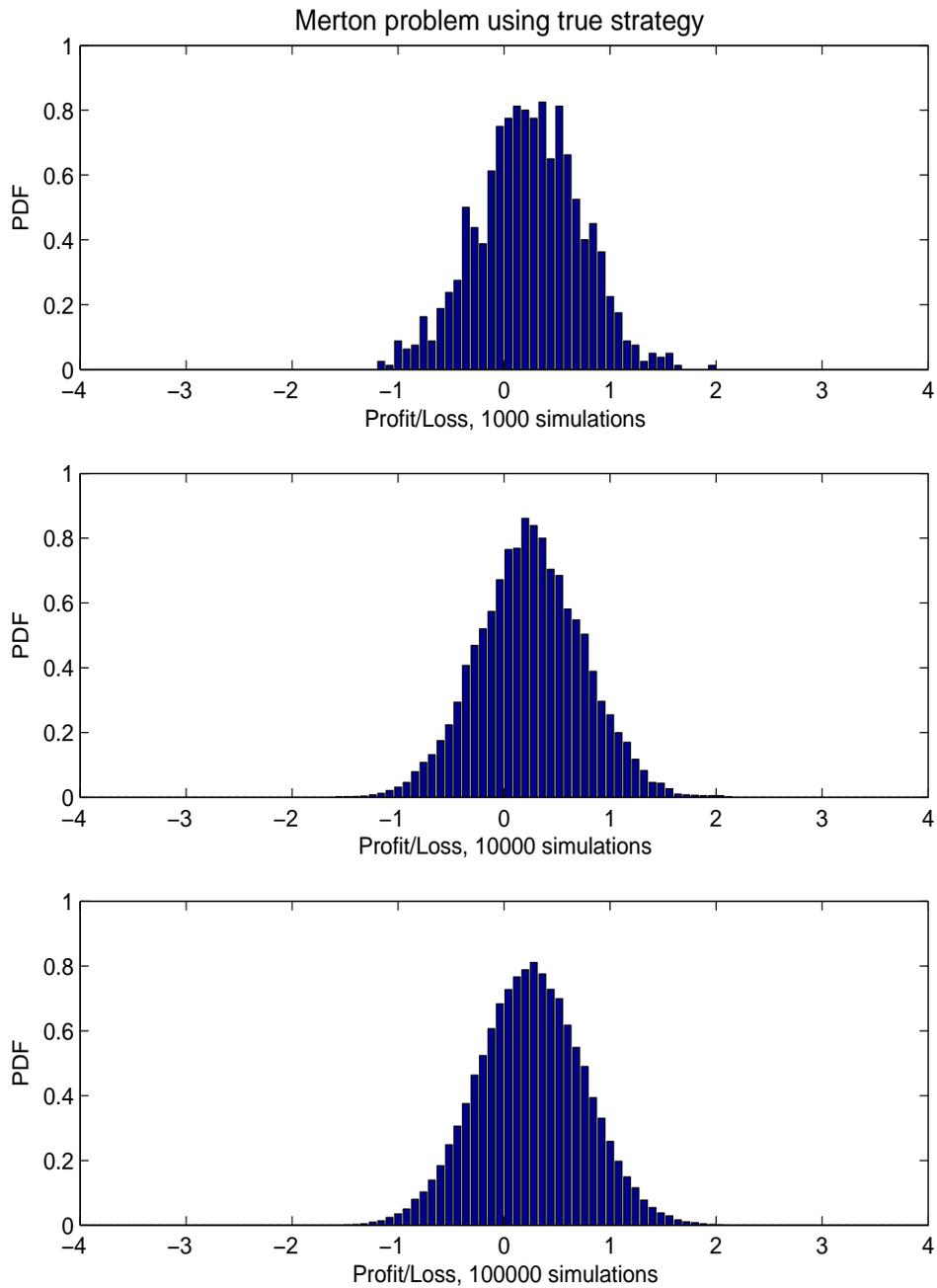}}
\caption{The profit/loss distributions of optimal investment portfolio, obtained as the exact solution
for Merton's problem with exponential utility, evaluated on
simulated stock prices following a geometric Brownian motion.}
\label{mertontruefig}
\end{figure}

\newpage
\begin{figure}[ht]
\epsfxsize=360pt
\epsfysize=500pt
\centerline{\epsffile{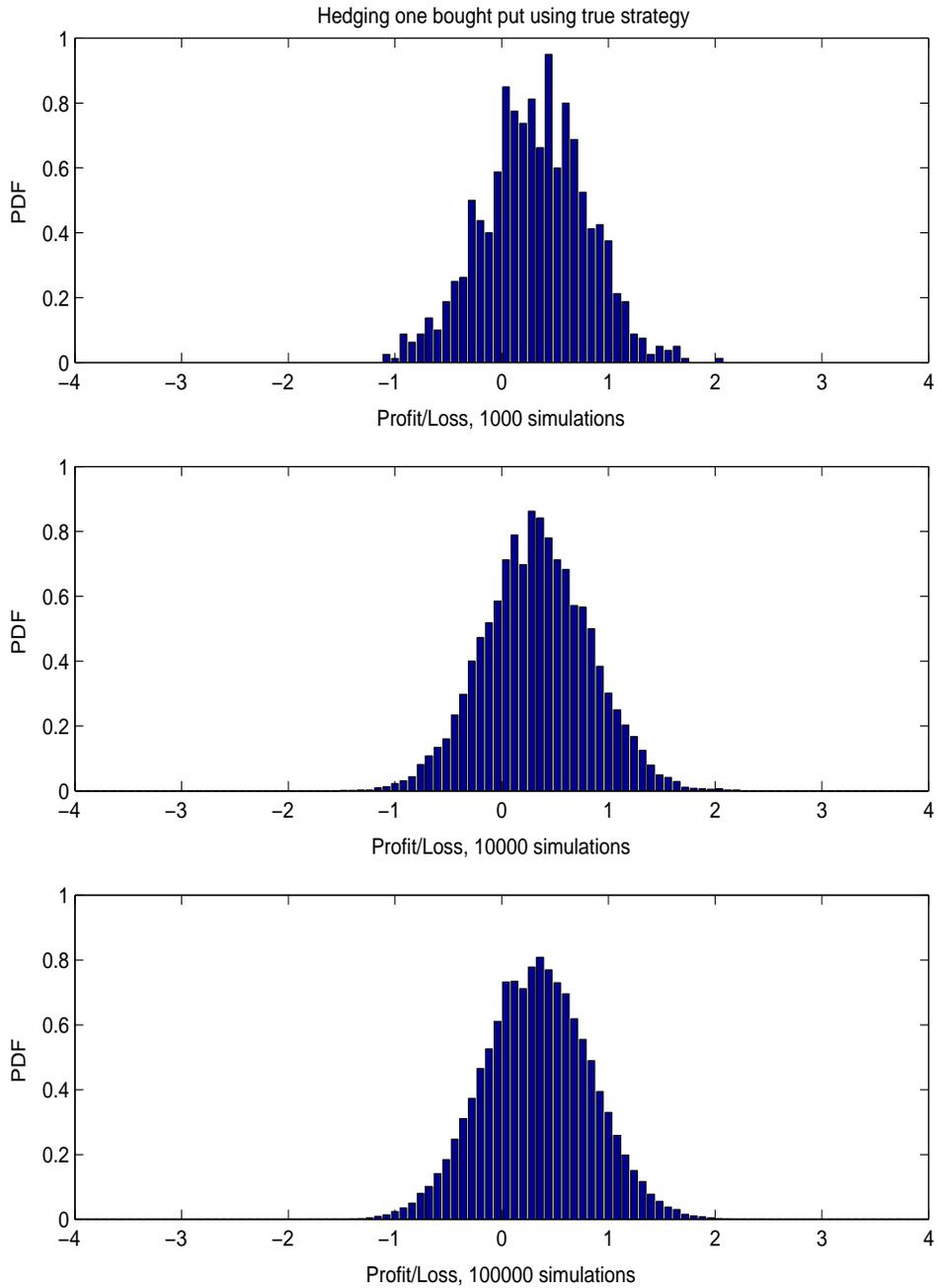}}
\caption{The profit/loss distributions of the optimal hedging portfolio for the buyer of one put option,
obtained from Black-Scholes delta hedging combined with Merton's problem with exponential utility, 
evaluated on simulated stock prices following a geometric Brownian
motion.}
\label{puttruefig}
\end{figure}

\newpage
\begin{figure}[ht]
\epsfxsize=380pt
\epsfysize=400pt
\centerline{\epsffile{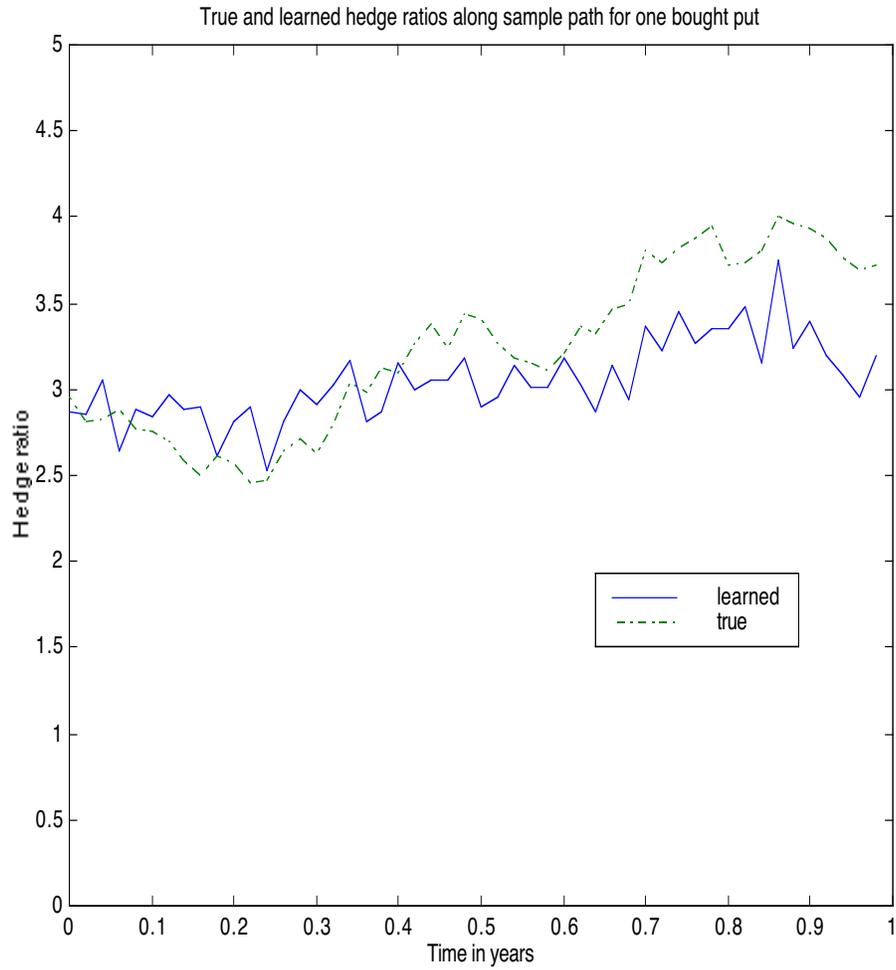}}
\caption{The hedge ratio (number of shares held) for the buyer of one put option on a simulated
sample path of duration one year, for which the option matures in--the--money. The solid line shows the 
strategy learned with N=100000; the broken line shows the theoretical Black-Scholes-Merton strategy.}
\label{hedgeplot}
\end{figure}

\end{document}